\documentclass{article}%

\usepackage[a4paper,left=3cm,right=2.5cm,top=3.0cm,bottom=3.0cm, headsep=0.5cm]{geometry}

\usepackage{amsmath}%
\usepackage{amsfonts}%
\usepackage{amssymb}%
\usepackage{amsthm}
\usepackage{graphicx}

\usepackage{xypic}
\input xy
\xyoption{all}

\newcommand{\N}{\mathbb{N}}

\newcommand{\R}{\mathbb{R}}

\newcommand{\D}{\mathrm{d}}

\newcommand{\M}{\mathbb{M}}

\newcommand{\conv}{\mathrm{conv\:}}

\newcommand{\bs}[1]{\boldsymbol{#1}}

\newenvironment{multi}{\left\{ \begin{array}{rcl}}{\end{array} \right.}

\newtheorem{theorem}{Theorem}[section]

\newtheorem{corollary}[theorem]{Corollary}

\newtheorem{definition}{Definition}[section]

\newtheorem{proposition}[theorem]{Proposition}

\begin{document}

\title{Convergence, Strong Law of Large Numbers, and Measurement Theory in the Language of Fuzzy Variables}
\author{\\ Adam Bzowski \thanks{e-mail: adamb@student.if.pw.edu.pl}\ , Micha\l  \ K. Urba\'{n}ski \thanks{e-mail: murba@if.pw.edu.pl} \\ \\
\normalsize Faculty of Physics, Warsaw University of Technology,\\ 
\normalsize ul. Koszykowa 75, 00-662 Warsaw, Poland}
\date{}
\maketitle

\begin{abstract}
In the paper we define the convergence of compact fuzzy sets as a convergence of $\alpha$-cuts in the topology of compact subsets of a metric space. Furthermore we define typical convergences of fuzzy variables and show relations with convergence of their fuzzy distributions. In this context we prove a general formulation of the Strong Law of Large Numbers for fuzzy sets and fuzzy variables with Archimedean $t$-norms. Next we dispute a structure of fuzzy logics and postulate a new definition of necessity measures. Finally, we prove fuzzy version of the Glivenko--Cantelli theorem and use it for a construction of a complete fuzzy measurement theory.
\end{abstract}

\section{Introduction}

\subsection{Motivations}

The application of fuzzy set theory in the description of empirical data becomes more popular in the recent years. A consistent theory requires algorithms which allow to construct a membership function of a given process with use of empirical data. This problem is completely solved in the theory of probability where by the Strong Law of Large Numbers a probability has an interpretation of frequency and a convergence of process estimators is guaranteed by Central Limit Theorem and Glivenko--Cantelli theorem. These theorems constitute a basis for further statistical analysis, construction of estimators and hypothesis verification.

The idea of fuzzy variables, counterparts of probabilistic random variables, reaches papers of Nahmias \cite{variables} and Rao and Rashed \cite{rao}. However, analogous theorems for fuzzy sets or fuzzy variables are not general enough. The most important and applicable result on the Law of Large Numbers is due to Triesch \cite{triesch}. However, this result was never interpreted in the language of fuzzy variables.

Before laws of large numbers may be proved, different types of convergence of fuzzy sets and fuzzy variables must be discussed. This requires to define at least a topological space of all fuzzy sets and fuzzy variables. Some general results on convergence of fuzzy sets are due to Diamond and Kloeden \cite{fuzzymetric} and Kaleva \cite{kaleva}. A convergence of fuzzy variables was used implicitely in the papers of Full\'{e}r \cite{fuller92}, Triesch \cite{triesch} and Williamson \cite{williamson}, but a general definition and dependencies between different convergences were not studied in detail.

Laws of Large Numbers were introduced in a few different versions. In the papers of Hong \cite{honglog}, Hong and Ro \cite{hong} and Triesch \cite{triescha} a pointwise convergence of fuzzy intervals is discussed. These results are based on analytical expressions for sums of fuzzy intervals developed by Full\'{e}r \cite{fuller92a}, Markov\'{a} \cite{markova}, Mesiar \cite{mesiar96}, Hong \cite{hongnote}, Hong and Hwang \cite{honghwang} and other which are valid only in special cases for example when their membership functions of the slopes are logarithmically concave. On the other hand a method of building up a membership function from a sample based on probability--possibility transformations discussed in the paper of Dubois \textit{et al.} \cite{dubois04} frequently leads to convex membership functions. 

A different approach to the Law of Large Numbers is presented in the papers of Full\'{e}r \cite{fuller92} and Triesch \cite{triesch}. We will show that a convergence analyzed in these papers is in fact a convergence of fuzzy variables. We will give an alternative proof of the theorem 3 in \cite{triesch}.

A discussion on convergence leads to a fuzzy counterpart of Glivenko--Cantelli theorem which we will proof in the last section. This theorem, crucial for a fuzzy data analysis, was never proved before.

In this paper we construct a complete theory of fuzzy variables. We discuss convergence of fuzzy variables and give a proof of the Strong Law of Large Numbers. Next, we construct an estimated membership function of a given process from a series of empirical data similarly to the procedure well-known from probability theory. Our fuzzy sets represent the whole process rather than a simple measurement result and describes both components of uncertainty, a random and a systematic one. This procedure does not require probability--possibility transformations or any probabilistic interpretation of fuzzy data.

Our philosophy is to look at fuzzy sets as a collection of $\alpha$-cuts rather than a membership function. Hence we would like to study a convergence of fuzzy sets as a convergence of their $\alpha$-cuts. This requires to define a topology on the space of $\alpha$-cuts which can be naturally done in case of compact fuzzy sets. In this fashion we want to estimate $\alpha$-cuts directly from empirical data rather than a membership function of a given fuzzy process and show a convergence of such a procedure to a fuzzy interval representing this process. In such an approach we use positional statistics such as modal value and median which are more natural in this context.

Finally, we must underline that our fuzzy variables are \textit{not} fuzzy random variables. The results for fuzzy random variables are similar to those known from probability and were intensively studied. Laws of Large Numbers are well-known in this case, e.g. \cite{teran}, \cite{klementlimit}.

\subsection{Definitions}

We use the following notations. By $\R$ we denote real numbers, $\N = \{ 0, 1, 2, \ldots \}$, $\left\lceil - \right\rceil$ and $\left\lfloor - \right\rfloor$ stand for floor and ceil functions respectively. For any set $\Omega$ by $2^\Omega$ we denote a set of all subsets of $\Omega$, for $A \subseteq \Omega$ we denote $A' = \Omega \backslash A$.

A normalized fuzzy set in the set $\Omega$ is a function $A : \Omega \rightarrow [0,1]$ such that there exists $x_0 \in \Omega$ with $A(x_0) = 1$. If $A(x) < 1$ for all $x \in \Omega$ then we call such a function a degenerated fuzzy set. In this paper by fuzzy set we always mean a normalized fuzzy set. For any $\alpha \in (0, 1]$ an $\alpha$-cut of a fuzzy set $A$ is a set $A^\alpha = \{ x \in \Omega : A(x) \geq \alpha \}$.


We will use following notations and definitions.
\begin{itemize}
	\item $\mathcal{FS}(\Omega)$ -- Set of all fuzzy sets in a set $\Omega$.
	\item $\mathcal{FK}(\Omega)$ -- Set of all compact fuzzy sets in a topological space $\Omega$, i.e. such fuzzy sets that for every $\alpha \in (0, 1]$ a set $A^\alpha$ is compact and nonempty.
	\item $\mathcal{FI}(\Omega)$ -- Set of all fuzzy intervals in a topological and linear space $\Omega$. $A \in \mathcal{FI}(\Omega)$ if and only if $A \in \mathcal{FK}(\Omega)$ and for every $\alpha \in (0, 1]$ a set $A^\alpha$ is convex.
	\item $\mathcal{FN}(\Omega)$ -- Set of all fuzzy numbers. $A \in \mathcal{FN}(\Omega)$ if $A$ is a fuzzy interval and there exists exactly one $x_0 \in \Omega$ such that $A(x_0) = 1$.
	\item $\mathcal{K}(\Omega)$ -- Set of all compact subsets of a metric space $\Omega$.
	\item $\mathcal{I} = \mathcal{KI}(\R)$ -- Set of all intervals.
\end{itemize}

For $x \in \mathcal{I}$ we denote $x = [\underline{x}, \overline{x}]$. Moreover for $A \in \mathcal{FI}(\R)$ we define $\overline{A}$ as a right slope of $A$, i.e. $\overline{A} = A|_{[\overline{A^1}, \infty)}$. Similarly for a left slope of $A$.

Every fuzzy set $A$ in a set $\Omega$ defines a fuzzy measure $\Pi_A$ on $\Omega$ defined for any $B \subseteq \Omega$ by $\Pi_A(B) = \sup A(B)$. Any measure obtained in this fashion is called a possibility measure or just a possibility. Any such a measure might be defined on a ring of sets $\mathcal{F} = 2^\Omega$ containing all subsets of $\Omega$, so we may not care about its domain.

Let $f : \Omega \rightarrow E$ be a function and $A \in \mathcal{FS}(\Omega)$. Then $f$ induces a map $f_\ast : \mathcal{FS}(\Omega) \rightarrow \mathcal{FS}(E)$ as $f_\ast(A)(y) = \sup A(f^{-1}(y))$. Similarly for any measure $\Pi_A$ on $\Omega$ the function $f$ induces a measure on $E$ as $f_\ast(\Pi_A)(B) = \Pi_A(f^{-1}(B))$.

Let $(\Omega, \Pi)$ be a set $\Omega$ with a possibility measure $\Pi$. A fuzzy variable on $\Omega$ with values in a set $E$ is any function $X : \Omega \rightarrow E$. Every fuzzy variable induces from $\Pi$ a possibility measure on $E$ called a distribution of $X$. A fuzzy set corresponding to this measure will be called a membership function of $X$. Fuzzy variables having fuzzy intervals as their membership functions will be called fuzzy interval variables while those with distributions being fuzzy numbers will be called fuzzy number variables.

A $t$-norm is a function $T : [0,1]^2 \rightarrow [0,1]$ which is commutative, associative, non-decreasing in both variables with $1$ being its neutral element. By its properties any $t$-norm may be extended to a function of countable number of arguments. $T$-norms generalize logical multiplication and serve as a definition of independence. We say that $B, C \subseteq \Omega$ are $T$-independent for some measure $\Pi_A$ if $\Pi_A(B \cap C) = T(\Pi_A(B), \Pi_A(C))$. A generalization to any number of events is routine. Fuzzy variables are $T$-independent if for all $B, C \subseteq \Omega$ sets $X^{-1}(B)$ and $X^{-1}(C)$ are $T$-independent with obvious extension for any number of fuzzy variables.

\section{Convergence}

\subsection{Convergence of fuzzy sets}

Let $(E, d_E)$ be a metric space and denote by $\mathcal{K}(E)$ a set of all compact subsets of $E$. The space $\mathcal{K}(E)$ is a metric space with a metric $d_{\mathcal{K}(E)}$ defined as
\begin{equation*}
d_{\mathcal{K}(E)}(K, L) = \max( \sup_{x \in K} \inf_{y \in L} d_X(x, y), \sup_{x \in L} \inf_{y \in K} d_E(x, y))
\end{equation*}

We define a convergence of general compact fuzzy sets as follows.

\begin{definition} \label{def:zbieg}
We say that a sequence $(A_n)_{n = 1}^\infty$ of compact fuzzy sets in a metric space $E$ converges to a compact fuzzy set $A$ levelwise if and only if for every $\alpha \in (0, 1]$ we have $A_n^\alpha \rightarrow A^\alpha$ in $\mathcal{K}(E)$.
\end{definition}

In the literature there are used also different types of convergence. Two most important are a pointwise convergence of membership functions, e.g. \cite{hong} and uniform convergence given by a metric $d(A, B) = \sup_{\alpha \in (0, 1]} d_{\mathcal{K}(E)}(A^\alpha, B^\alpha)$, e.g. \cite{kaleva}. However, the first one does not agree with philosophy of fuzzy sets seeing as a collection of $\alpha$-cuts, while the latter is too strong if we do not assume compactness of supports. Our definition \ref{def:zbieg} will occur to be more natural for fuzzy sets and consistent with fuzzy variables convergence. Note, that in this topology maps $A \mapsto A^\alpha$ are continuous for all $\alpha \in (0, 1]$. The following proposition links our definition of convergence of compact fuzzy sets with pointwise convergence of their membership functions in some special case.

\begin{proposition} \label{zbieg}
Let $(A_n)_{n=1}^\infty$ be a sequence of compact fuzzy sets in a complete metric space $E$.
\begin{enumerate}
	\item If for all $\alpha \in (0, 1]$, $A_{n+1}^\alpha \subseteq A_n^\alpha$ for all $n$, then $(A_n)_{n=1}^{\infty}$ converges to a compact fuzzy set $A$ with $\alpha$-cuts given by $A^\alpha = \bigcap_{n=1}^\infty A_n^\alpha$.
	\item If for all $\alpha \in (0, 1]$, $A_{n+1}^\alpha \subseteq A_n^\alpha$ for all $n$, then for every $x \in E$ we have $A_n(x) \rightarrow A(x)$ with $A$ obtained from the previous point.
	\item If for every $x \in E$ a sequence $(A_n(x))_{n=1}^\infty$ is non-increasing and converges to a membership function of a compact fuzzy set $A$, then $A_n \rightarrow A$ in the sense of the definition \ref{def:zbieg}.
 \end{enumerate}
\end{proposition}

\begin{proof}
\begin{enumerate}
	\item By definition we must show that if $(K_n)_{n=1}^\infty$ is a sequence of compact subsets of $E$ such that $K_{n+1} \subseteq K_n$ for all $n$, then $(K_n)_{n=1}^\infty$ converges to $K = \bigcap_{n=1}^\infty K_n$, which is non-empty by completeness. Observe that for a given $N$ we have $d(K, K_N) = \sup_{x \in K_N} d(K, x)$. Now assume that $(K_n)_{n=1}^\infty$ does not converge to $K$, i.e. there exists $\epsilon > 0$ such that for all $N$ there exists $n \geq N$ and $x \in K_n$ with $d(x, K) > \epsilon$. Firstly $x \notin K$, but on the other hand $x \in K_n$ for all $n$. Contradiction.
	\item By definition $A_n(x) = \sup \{ \alpha : x \in A_n^\alpha \}$ and similarly for $A(x)$. Then, assuming $\sup \emptyset = 0$ we have
\begin{equation*}
A(x) = \sup \{ \alpha : x \in \bigcap_{n=1}^\infty A_n^\alpha \} = \sup \bigcap_{n=1}^\infty \{ \alpha : x \in A_n^\alpha \} = \lim_{n \rightarrow \infty} \sup \{ \alpha : x \in A_n^\alpha \} = \lim_{n \rightarrow \infty} A_n(x)
\end{equation*}
	The limit exists since $A_n(x)$ is non-increasing.
\item It is clear that $A^\alpha_{n+1} \subseteq A^\alpha_n$ for all $\alpha \in (0, 1]$ and hence by point 1 we have $A_n \rightarrow A$.
\end{enumerate}

\end{proof}

\subsection{Convergence of fuzzy variables}

Firstly, let us recall a definition of fuzzy variable and its membership function.
\begin{definition}
Let $(\Omega, \Pi)$ be a set with a possibility measure given by a fuzzy set $A$ in $\Omega$. A fuzzy variable $X$ with values in a set $E$ is any function $X : \Omega \rightarrow E$. A possibility measure $X_\ast \Pi$ on $E$ induced by $X$ is called a distribution of $X$ while a fuzzy set $X_\ast A$ in $E$ determining the distribution is called a membership function of $X$.
\end{definition}

\begin{definition}
Let $(X_n)_{n=1}^\infty$ be a sequence of fuzzy variables defined on a space $(\Omega, \Pi)$ with values in $\R^n$ (in general in some Polish space). We define following types of convergence.
\begin{enumerate}
	\item Almost sure convergence. $\Pi( \{ \omega \in \Omega : \lim_{n \rightarrow \infty} X_n(\omega) \neq X(\omega) \} ) = 0$,
	\item Weak almost sure convergence. $\Pi( \{ \omega \in \Omega : \lim_{n \rightarrow \infty} X_n(\omega) = X(\omega) \} ) = 1$,
	\item Convergence in measure. For every $\epsilon > 0$ we have $\lim_{n \rightarrow \infty} \Pi( |X_n - X| \geq \epsilon) = 0$.
	\item Weak convergence in measure. For every $\epsilon > 0$ we have $\lim_{n \rightarrow \infty} \Pi( |X_n - X| \leq \epsilon) = 1$.
\end{enumerate}
Moreover if all $X_n$ are compact fuzzy variables then we say that $(X_n)_{n=1}^\infty$ converges to a compact fuzzy variable $X$ in distribution if and only if membership functions of $X_n$ converge to the a membership function of $X$ in the sense of definition \ref{def:zbieg}.
\end{definition}

Observe that a convergence in measure is equivalent to the Full\'{e}r's definition \cite{fuller92} when a limit is one-point fuzzy variable. We would prefer not to use necessity measures now because we are going to discuss them later. Note that a convergence in distribution is similar to a definition given by Williamson \cite{williamson} with pointwise convergence exchanged by levelwise.

\begin{proposition} \label{zbiegX}
Let $(X_n)_{n=1}^\infty$ be a sequence of fuzzy variables with values in $\R^n$.
\begin{enumerate}
	\item[a)] There are following dependencies between particular types of convergence: convergence in measure $\Rightarrow$ almost sure convergence $\Rightarrow$ weak almost sure convergence $\Rightarrow$ weak convergence in measure.
	\item[b)] If $c$ is a "non-fuzzy" one-point fuzzy variable, then convergence of distributions is equivalent with convergence in measure.
\end{enumerate}
\end{proposition}

\begin{proof}
For the first arrow we have the following reasoning. Assume for every $\epsilon > 0$ we have $\lim_{n \rightarrow \infty} \Pi( |X_n - X| \geq \epsilon) = 0$. By definition $(X_n)_{n=1}^\infty$ converges almost surely if for any $\epsilon > 0$ we have 
\begin{equation*}
\Pi \left( \bigcap_{N = 1}^\infty \bigcup_{n=N}^\infty \{ \omega \in \Omega : |X_n(\omega) - X(\omega)| \geq \epsilon \} \right) = 0
\end{equation*}
Assume $(B_n)_{n=1}^\infty$ a sequence of subsets of $\Omega$ such that $B_{n+1} \subseteq B_n$. Then a sequence $(\Pi(B_n))_{n=1}^\infty$ is non-increasing and has a limit. Moreover $\Pi \left( \bigcap_{n=1}^\infty B_n \right) \leq \lim_{n \rightarrow \infty} \Pi(B_n)$. Hence
\begin{eqnarray*}
&& \Pi \left( \bigcap_{N = 1}^\infty \bigcup_{n=N}^\infty \{ |X_n - X| \geq \epsilon \} \right)\\
& \leq & \lim_{N \rightarrow \infty} \Pi \left(\bigcup_{n=N}^\infty \{ |X_n - X| \geq \epsilon \} \right) \\
& = & \lim_{N \rightarrow \infty} \sup_{n \geq N} \Pi(|X_n - X| \geq \epsilon) = 0 \\
\end{eqnarray*}
since by assumption the limit exists.

The second arrow is obvious. For the third observe that for every $\epsilon > 0$ we have
\begin{eqnarray*}
1 & = & \Pi \left( \bigcup_{N = 1}^\infty \bigcap_{n=N}^\infty \{ |X_n - X| \leq \epsilon \} \right) \\
& = & \sup_{N} \Pi \left(\bigcap_{n=N}^\infty \{ |X_n - X| \leq \epsilon \} \right) \\
& \leq & \sup_{N} \lim_{n \geq N} \Pi(|X_n - X| \leq \epsilon) \\
& = & \lim_{n \rightarrow \infty} \Pi(|X_n - X| \leq \epsilon)
\end{eqnarray*}

For b) part denote by $A_{X_n}$ membership functions of $X_n$ and observe the following
\begin{equation*}
\Pi( |X_n - c| \geq \epsilon) \leq \delta \; \iff \; A_{X_n}^\delta \subseteq B(c,\epsilon)
\end{equation*}
where $B(c,\epsilon)=\{x\in\R^n: |x-c| \leq \epsilon\}$ from which thesis follows.

\end{proof}

\section{Tools}

\subsection{$T$-norms, conorms and measures}

The following well known facts can be found in Alsina \textit{et al.} \cite{alsina}, Klement \cite{klement} or Acz\'{e}l \cite{aczel}.

$T$-norm $T$ is called Archimedean if for any $x, y \in (0, 1)$ there exists $n \in \N$ such that $T(\underbrace{x, \ldots, x}_{n}) < y$. Every continuous Archimedean $t$-norm has a representation $T(x, y) = g^{[-1]}(g(x) + g(y))$ where $g : [0,1] \rightarrow [0, \infty]$ is a strictly decreasing and convex function such that $g(1) = 0$ called an additive generator of $T$. Here $g^{[-1]} : [0, \infty] \rightarrow [0,1]$ is a pseudo-inverse, i.e. $g^{[-1]}(y) = g^{-1}(y)$ for $y \in [0, g(0)]$ and $0$ otherwise. This representation is unique up to a multiplication by a positive constant. A $t$-norm $T$ is called strict if it is strictly increasing. A continuous Archimedean $t$-norm is strict if and only if $g(0) = +\infty$. 

If $g$ is an additive generator of a continuous Archemedean $t$-norm $T$, than $h(x) = e^{-g(x)}$ is a multiplicative generator of $T$, i.e. $T(x,y) = h^{[-1]}(h(x)h(y))$. Here $h : [0,1] \rightarrow [0, 1]$ and a psuedo-inverse $h^{[-1]} : [0, 1] \rightarrow [0,1]$ is defined as $h^{[-1]}(y) = h^{-1}(y)$ for $y \in [h(0), 1]$ and $0$ otherwise. A continuous $t$-norm $T$ is Archimedean if and only if it admits a representation by a multiplicative generator $h : [0,1] \rightarrow [0,1]$. If $h_1$ and $h_2$ are two generators of $T$ then $h_2(x) = (h_1(x))^\alpha$ and $h_1^{[-1]}(y) = h_2^{[-1]}(x^\alpha)$ for some $\alpha > 0$. Hence a multiplicative generator is of the form $h(x) = e^{-\alpha g(x)}$ and so $h(1) = 1$. $T$ is strict if and only if $h(0) = 0$.

A $t$-conorm $S$, also known as $s$-norm, is a function $S : [0,1]^2 \rightarrow [0,1]$ which is commutative, associative, non-decreasing in both variables with $0$ being its neutral element. By its properties any $t$-conorm may be extended to a function of countable number of arguments. $T$-conorms generalize logical addition and in the theory of decomposable measures serve as a definition of generalized addition. A measure $\mu : \mathcal{F} \rightarrow [0,1]$ is called $S$-decomposable and normalized if $\mu(\emptyset) = 0$, $\mu(X) = 1$ and for any $B, C \in \mathcal{F}$ such that $B \cap C = \emptyset$, $\mu(B \cup C) = S(\mu(B), \mu(C))$. A generalization to the countable number of sets follows from properties of $t$-conorms. Further discussion on decomposable measures can be found in \cite{null} and \cite{decomp}

A continuous $t$-conorm $S$ is called Archimedean if $T(x, y) = 1 - S(1 - x, 1 - y)$ is an Archimedean $t$-norm. Every continuous Archimedean $t$-conorm has a representation $S(x, y) = g^{[-1]}(g(x) + g(y))$ where $g : [0,1] \rightarrow [0, \infty]$ is an increasing function such that $g(0) = 0$ called an additive generator of $S$. This representation is unique up to a multiplication by a positive constant. A $t$-conorm $S$ is strict if it is strictly increasing. A continuous Archimedean $t$-conorm is strict if and only if $g(1) = +\infty$.

Let $A \in \mathcal{FS}(X)$ and $B \in \mathcal{FS}(Y)$. A $T$-product of fuzzy sets $A \otimes B \in \mathcal{FS}(X \times Y)$ for some $t$-norm $T$ is a fuzzy set defined as $(A \otimes B)(x, y) = T(A(x), B(y))$. A fuzzy measure coming from $A \otimes B$ is a unique possibility product measure.

\subsection{Nguyen-Full\'{e}r-Keresztfalvi theorem}

If $X, Y$ are two $T$-independent fuzzy variables with values in $E$ and memberrship functions $A_X$ and $A_Y$ respectively and $f : E \times E \rightarrow E'$ is a function, then a fuzzy variable $Z = f(X,Y)$ has a membership function $A_Z$ given by Zadeh Extension Principle
\begin{equation*}
A_Z(z) = \sup \{ T(A_X(x), A_Y(y)) : f(x,y) = z \}
\end{equation*}
By Nguyen-Full\'{e}r-Keresztfalvi theorem (NFK theorem) \cite{fuller91} if $f$ is continuous, $A_X$, $A_Y$, $T$ are upper semicontinuous and $A_X$, $A_Y$ are compactly supported then
\begin{equation*}
A_Z^\alpha = \bigcup_{T(\xi, \eta) \geq \alpha} f(A_X^\xi, A_Y^\eta)
\end{equation*}
with $f$ applied in a set theoretic sense. This theorem can be generalized to the following version

\begin{theorem}[Nguyen-Full\'{e}r-Keresztfalvi theorem] \label{nfk}
Let $X, Y, Z$ be topological spaces and let in $Z$ every one-point subset be closed. Let $f : X \times Y \rightarrow Z$ be a continuous map, $T$ an upper semicontinuous $t$-norm and $A \in \mathcal{FK}(X)$, $B \in \mathcal{FK}(Y)$. Then
\begin{equation*}
[f_\ast(A, B)]^\alpha = \bigcup_{T(\xi, \eta) \geq \alpha} f(A^\xi, B^\eta)
\end{equation*}
for all $\alpha \in (0, 1]$ where $f_\ast : \mathcal{FK}(X) \times \mathcal{FK}(Y) \rightarrow \mathcal{FK}(Z)$ is induced from $f$ by Zadeh Extension Principle.
\end{theorem}
\begin{proof}
We need to show three points.
\begin{enumerate}
	\item $A$ and $B$ have upper semicontinuous membership functions. Indeed, take a sequence $(x_n)_{n=1}^\infty$ of points of $X$ such that $x_n \rightarrow x$ and $(A(x_n))_{n=1}^\infty$ is non-increasing. Then if for some $N \in \N$ and $\alpha \in (0, 1]$ we have $x_N \in A^\alpha$, then $x_n \in A^\alpha$ for all $n \leq N$ and hence $A(x) \leq \lim_{n \rightarrow \infty} A(x_n)$. On the other hand by compactness if all $x_n \in A^\alpha$ then $x \in A^\alpha$, so $\lim_{n \rightarrow \infty} A(x_n) \leq A(x)$.
	\item By theorem 1 of \cite{fuller91} we need to show that a function $\varphi(x, y) = (A \otimes B)(x, y)$ attains its supremum on $f^{-1}(z)$ for every $z \in Z$. If this supremum is non-zero then there exists $(x_0, y_0) \in f^{-1}(z)$ such that $\varphi(x_0, y_0) = \epsilon > 0$. Because for $T(\xi, \eta) \geq \epsilon$ it is necessary that $\xi \geq \epsilon$ and $\eta \geq \epsilon$ we see that the supremum is attained in a set $A^\epsilon \times B^\epsilon$ which is compact. By assumptions $A^\epsilon \times B^\epsilon \cap f^{-1}(z)$ is compact and hence $\varphi$ attains its maximum. 
	\item The image of $f_\ast$ lies in $\mathcal{FK}(Z)$. Because by previous point $\varphi$ attains its supremum, then we have by definition
\begin{equation*}
[f(A, B)]^\alpha = \{ z \in Z : \sup_{f(x,y) = z} \varphi(x, y) \geq \alpha \} = f( \varphi^{-1}([\alpha, 1]))
\end{equation*}
	A set $\varphi^{-1}([\alpha, 1])$ is a closed subset of a compact set $A^\alpha \times B^\alpha$ and hence $[f(A, B)]^\alpha$ is compact.
\end{enumerate}
\end{proof}

In the formulation of Full\'{e}r and Ketereszfalvi the assumption on $Z$ is missing. Now we have a few corollaries which are very important in further analysis.

\begin{corollary}
Let $X$ and $Z$ be topological spaces and let in $Z$ every one-point set be closed. Let $f : X \rightarrow Z$ be a continuous function and $A \in \mathcal{FK}(X)$. Then for an induced map $f_\ast : \mathcal{FK}(X) \rightarrow \mathcal{FK}(Z)$ we have $(f_\ast A)^\alpha = f(A^\alpha)$.
\end{corollary}
\begin{proof}
Take for $Y$ a one-point topological space and as $B$ a one-point fuzzy set in $Y$ and apply theorem \ref{nfk}.
\end{proof}

\begin{corollary}
If $T_1$ and $T_2$ are two upper semicontinuous $t$-norms such that $T_1 \leq T_2$ and $X$, $Y$, $Z$, $A$, $B$, $f$ as in the NFK theorem \ref{nfk} then $f_\ast^{(T_1)}(A, B) \subseteq f_\ast^{(T_2)}(A, B)$ where $f_\ast^{(T)}$ stands for a function induced from $f$ by means of Zadeh Extension Principle with use of a $t$-norm $T$.
\end{corollary}
\begin{proof}
If $T_1 \leq T_2$ then $\{ (\xi, \eta) : T_1(\xi, \eta) \geq \alpha \} \subseteq \{ (\xi, \eta) : T_2(\xi, \eta) \geq \alpha \}$ and hence by monotonicity of $t$-norms $[f_\ast^{(T_1)}(A, B)]^\alpha \subseteq [f_\ast^{(T_2)}(A, B)]^\alpha$.
\end{proof}

\begin{corollary}
Let $T$ be an Archimedean $t$-norm with an additive generator $g$ and let $X$, $Y$, $Z$, $A$, $B$, $f$ as in the NFK theorem \ref{nfk}. Then
\begin{equation*}
[f_\ast(A, B)]^\alpha = \bigcup_{g(\xi) + g(\eta) \leq g(\alpha)} f(A^\xi, B^\eta)
\end{equation*}\end{corollary}
\begin{proof}
$g \circ g^{[-1]}$ is not an identity map. If $(\xi, \eta)$ are such that $g(\xi) + g(\eta) \leq g(0)$ then $g(\xi) + g(\eta) = g(\alpha) \iff T(\xi, \eta) = \alpha$. If $g(\xi) + g(\eta) > g(0)$ then $T(\xi, \eta) = 0$ but the equation $g(\xi) + g(\eta) \leq g(\alpha)$ does not have a solution. Hence the only difference is for $0$-cut and so is meaningless.
\end{proof}

This corollary can be easily extended to a formula valid for any finite number of fuzzy sets.

\subsection{Fuzzy modal value}

A fuzzy modal value is a counterpart of a probabilistic expectation value, however, unlike in probability theory, it is a set rather than a number. It is defined as follows.

\begin{definition}
A fuzzy modal map $\M : \mathcal{FS}(\Omega) \rightarrow 2^\Omega \backslash \{ \emptyset \} \subseteq \mathcal{FS}(\Omega)$ is a function defined as $\M A = A^1$. Similarly for fuzzy variables we define $\M$ as the operator with values in $2^\Omega \backslash \{ \emptyset \}$ as $\M X = \M A_X$ where $A_X$ is a membership function of $X$.
\end{definition}

If a fuzzy set $A$ is a fuzzy number or a fuzzy variable $X$ is a fuzzy number variable, then we treat $\M A$ and $\M X$ as a point (number) rather than a one-point set.

\begin{proposition}
Let $f : E \rightarrow E'$ be a continuous function between two topological spaces with $E'$ such that all one-point sets are closed. Let $A, B \in \mathcal{FK}(E)$ and $X : \Omega \rightarrow E$ be a fuzzy variable with distribution given by $A$. We can treat $\M A$ as a fuzzy set and hence
\begin{enumerate}
  \item $\M (\mathcal{FK}(E)) = \mathcal{K}(E) \subseteq \mathcal{FK}(E)$
	\item $\M f_\ast A = f_\ast \M A$,
	\item $\M T(A, B) = T(\M(A), \M(B))$,
	\item $\M f(X) = f(\M X)$.
\end{enumerate}
\end{proposition}
\begin{proof}
By NFK theorem we have the following
\begin{enumerate}
  \item by definition,
	\item $\M f_\ast A = (f_\ast A)^1 = f(A^1) = f_\ast \M A$,
	\item $\M T(A, B) = T(A, B)^1 = (A \cap B)^1 = A^1 \cap B^1 = T(A^1, B^1) = T(\M A, \M B)$,
	\item $\M f(X) = A_{f(X)}^1 = (f_\ast A_X)^1 = f(A_X^1) = f(\M X)$.
\end{enumerate}
\end{proof}

Note that properties 2 and 4 are true for any continuous functions so are more general than their counterparts in probability theory because expectation value is a linear operator.

\subsection{Mean and median}

Let us define two most important functions of type $\R^n \rightarrow \R$ for any $n \in \N$ which will br used in the further analysis.
\begin{eqnarray*}
Ave(x_1, \ldots, x_n) & = & \frac{1}{n} \sum_{j=1}^n x_j \\
Med(x_1, \ldots, x_{2n+1}) & = & median(x_1, \ldots, x_{2n+1}) \\
\end{eqnarray*}

\begin{proposition} \label{dlaMC}
Let $T$ be a continuous Archimedean $t$-norm with an additive generator $g$. Let $\alpha \in (0, 1]$, $n \in \N$ and $A$ be a compact fuzzy set in $\R$. Denote $N = 2n + 1$. Then
\begin{equation*}
Med_\ast(\underbrace{A, \ldots, A}_{N})^\alpha = A^{g^{[-1]}\left( \frac{2}{N+1} g(\alpha) \right) }
\end{equation*}
\end{proposition}
\begin{proof}
The median induces by NFK theorem a map
\begin{equation*}
Med_\ast(\underbrace{A, \ldots, A}_{N})^\alpha = \bigcup_{g(\xi_1) + \ldots + g(\xi_N) \leq g(\alpha)} Med(A^{\xi_1}, \ldots, A^{\xi_N})
\end{equation*}
with $Med(A^{\xi_1}, \ldots, A^{\xi_N})$ being such a set among $A^\xi_j$ which contains and is contained in exactly $n$ different sets $A^{\xi_j}$.

Let $\xi_1, \ldots, \xi_N$ satisfy $g(\xi_1) + \ldots + g(\xi_N) \leq g(\alpha)$. For any given sets $A^{\xi_1}, \ldots, A^{\xi_N}$ we may order them and so assume $A^{\xi_1} \subseteq \ldots \subseteq A^{\xi_N}$. Now if we set $\tilde{\xi}_1 = \ldots = \tilde{\xi}_n = 1$ and $\tilde{\xi}_{n+1} = \xi_{n+1}, \ldots, \tilde{\xi}_{N} = \xi_{N}$ then $g(\tilde{\xi}_1) + \ldots + g(\tilde{\xi}_N) \leq g(\alpha)$ and $Med(A^{\xi_1}, \ldots, A^{\xi_N}) = Med(A^{\tilde{\xi}_1}, \ldots, A^{\tilde{\xi}_N})$. Because for every $\xi \leq \eta$ we have $\M A \subseteq A^\eta \subseteq A^\xi$, hence $Med(A^{\tilde{\xi}_1}, \ldots, A^{\tilde{\xi}_N})$ is such a set among $A^{\tilde{\xi}_{n+1}}, \ldots, A^{\tilde{\xi}_{N}}$ which is contained in all the others. The largest possible such a set requires all $A^{\tilde{\xi}_j}$ to be equal with minimal $\tilde{\xi}_j$ for $j \geq n+1$ and hence $\tilde{\xi}_j = g^{[-1]}\left( \frac{1}{n+1} g \left( \alpha \right) \right)$. We obtain thesis by substitution $n + 1 = \frac{N+1}{2}$.
\end{proof}

In the calculation of median we have assumed $N$ to be an odd number. Now we postulate a result for any $n \in \N$ as
\begin{equation*}
Med_\ast(\underbrace{A, \ldots, A}_{n})^\alpha = A^{g^{[-1]}\left( \frac{2}{n+1} g(\alpha) \right) }
\end{equation*}

\section{Laws of Large Numbers}

\begin{theorem}[Law of Large Numbers for median] \label{bLoLNmed}
Let $T$ be a continuous Archimedean $t$-norm with an additive generator $g$ and let $A \in \mathcal{FK}(\R)$ For any $n \in \N$ define compact fuzzy sets
\begin{equation*}
M_n = Med_\ast(\underbrace{A, \ldots, A}_{n})
\end{equation*}
Then $(M_n)_{n=1}^\infty$ converges to $\M A$.
\end{theorem}
\begin{proof}
For $M_n$ by proposition \ref{dlaMC} we have
\begin{equation*}
Med_\ast(\underbrace{A, \ldots, A}_{N})^\alpha = A^{g^{[-1]}\left( \frac{2}{N+1} g(\alpha) \right)}
\end{equation*}
It follows that $M_{n+1}^\alpha \subseteq M_n^\alpha$ and $\M M_n = \M A$ so by propositon \ref{zbieg} $M_n \rightarrow M$ with $M^\alpha = \bigcap_{n=1}^\infty M_n^\alpha \supseteq \M A$. On the other hand if there exists $x \in \R$ such that $x \in M_n^\alpha$ for all $n$ then $x \in A^\beta$ for all $\beta < 1$ and hence $x \in \M A$. Hence $(M_n)_{n=1}^\infty$ converges to $\M A$.
\end{proof}

\begin{theorem}[Law of Large Numbers for mean] \label{bLoLNmean}
Let $T$ be a continuous Archimedean $t$-norm with an additive generator $g$ and let $A \in \mathcal{FI}(\R)$ For any $n \in \N$ define fuzzy intervals
\begin{equation*}
A_n = Ave_\ast(\underbrace{A, \ldots, A}_{n}) \\
\end{equation*}
Then $(A_n)_{n=1}^\infty$ converges to $\M A$.
\end{theorem}
\begin{proof}
First assume $A$ is continuous and concave, i.e. $\underline{A}$ and $\overline{A}$ are continuous and concave functions. Then by means of works \cite{honghwank}, \cite{hongnote}, \cite{markova}, \cite{mesiar96} or \cite{triescha} we find
\begin{equation*}
Ave_\ast(\underbrace{A, \ldots, A}_{n})^\alpha = A^{g^{[-1]}\left( \frac{1}{n} g(\alpha) \right)}
\end{equation*}

Assume the right slope $\overline{A}$ of $A$ is not concave and $a = \overline{A^1}$. Denote $\phi(\alpha) = \overline{A^\alpha}$ and note that $\phi : (0,1] \rightarrow [a, \infty)$ is a left-continuous non-increasing function with $\phi(1) = a$. Now set $\alpha \in (0, 1]$ and denote $Y = (0,1] \times [a, \infty)$ and
\begin{equation*}
C = \{ (\beta, x) \in Y : x \leq \phi(\alpha) \}
\end{equation*}
and $\tilde{C} = \conv(C)$ a convex hull of $C$. For any $x_0 \in [a, \infty)$ denote $\tilde{C}_{x_0} = \{ \beta \in (0, 1] : (\beta, x_0) \in \tilde{C} \}$. Note that $C$ is closed in $Y$, so is $\tilde{C}$ and $\tilde{C}_{x_0}$ for any $x_0 \in [a, \infty)$ and hence denote $H(x) = \max \tilde{C}_{x}$. Moreover $H(x) > H(y)$ if $x < y$ and by definition if $(1, x) \in \tilde{C}$ then $x = a$ and hence for any $x \in (a, \infty)$ we have $H(x) < 1$. $H$ is continuous since any line passing through $(x, H(x))$ and $(\phi(\alpha), \alpha)$ meets $C_y$ in some point and
\begin{equation*}
|H(x) - H(y)| = |\max \tilde{C}_{x} - \max \tilde{C}_{y}| \leq \frac{H(\max(x, y)) - \alpha}{\phi(\alpha) - \max(x, y)} |x - y|
\end{equation*}
Now we can define a fuzzy interval $\tilde{H} : [a, \infty) \rightarrow [0,1]$ as follows
\begin{equation*}
\tilde{H}(x) = \left\{ \begin{array}{rcl} 
	H(x) & \mathrm{if} & x \leq \phi(\alpha) \\
	\overline{A}(x) & \mathrm{if} & x > \phi(\alpha) \end{array} \right.
\end{equation*}
$\tilde{H}$ defines a fuzzy interval on $[a, \infty)$ with continuous membership function which is concave, and in particular logarithmically concave on $[a, \phi(\alpha)]$.

Now carry out the same procedure for the second slope of $A$, concatenate the fuzzy sets in order to obtain a fuzzy interval $B$ on $\R$ with continuous membership function and concave on $A^\alpha$. Observe that $A \subseteq B$ in order to get for any $\alpha \in (0, 1]$ in a similar way as in the proof of theorem \ref{bLoLNmed}
\begin{equation*}
A_n^\alpha \subseteq B_n^{\alpha^{1/n}} \stackrel{n \rightarrow \infty}{\longrightarrow} \M B = \M A
\end{equation*}
and finish the proof.
\end{proof}

By proposition \ref{zbiegX} for the convergence of fuzzy number variables it is enough to prove the convergence of distributions. Hence we have

\begin{theorem}[Strong Law of Large Numbers] \label{SLoLN}
Let $T$ be a continuous Archimedean $t$-norm and let $A \in \mathcal{FN}(\R)$. Let $(X_n)_{n=1}^\infty$ be a sequence of $T$-independent fuzzy number variables defined on a space $(\Omega, \Pi)$ with distributions given by $A$. For any $n \in \N$ define fuzzy number variables
\begin{eqnarray*}
A_n & = & Ave(X_1, \ldots, X_n) \\
M_n & = & Med(X_1, \ldots, X_n)
\end{eqnarray*}
Then $(A_n)_{n=1}^\infty$ and $(M_n)_{n=1}^\infty$ converge in distribution, in measure, and almost surely to $\M A$ viewed as a one-point fuzzy variable.
\end{theorem}

This theorem was proved by Triesch in paper \cite{triesch}. However, it was neither interpreted in a context of fuzzy variables nor compared with dependencies of fuzzy sets.

\section{Set complements, necessity measures and logics}

Consider the following facts.
\begin{itemize}
	\item We may treat fuzzy sets as a dual theory to probability by means of the formula $P(A) = \int A(x) \D P(x)$ where $A$ is a fuzzy set and $P$ a probability measure. However, a definition of Lebesgue-like integral involves addition and multiplication, so a choice of a $t$-norm should define a choice of a $t$-conorm representing addition.
	\item For a given possibility measure $\Pi$ there exists a necessity measure $\Pi'$ given by a formula $\Pi'(A) = 1 - \Pi(A')$. However this definition should depend on a $t$-norm since for example a relation $A' \cap B = B \backslash A$ involves $t$-norm.
	\item In the sense of necessity it is reasonable to say that a possibility of \textit{not} $A$ is equal to $1 - \Pi(A)$ and probably not identical with $\Pi(A')$.
\end{itemize}

\subsection{Logics and norms}

We would like to define a logic-like structure of a type $([0,1], \oplus, \otimes, \neg)$ where $\oplus$, $\otimes$ and $\neg$ stand for logical \textit{or}, \textit{and} and negation. Moreover we would like to say that a measure $\mu$ has its values in this structure in such a sense that
\begin{enumerate}
	\item if $A \cap B = \emptyset$ then $\mu(A \cup B) = \mu(A) \oplus \mu(B)$, i.e., $\mu$ is $\oplus$-decomposable measure,
	\item if $\mu(A \cap B) = \mu(A) \otimes \mu(B)$ then we say that $A$ and $B$ are independent,
	\item $\mu(A') = \neg \mu(A)$, where $A' = X \backslash A$.
\end{enumerate}

Due to results of \cite{alsina} it is impossible for this structure to be a Boolean logic, since if $\oplus$ and $\otimes$ are mutually distributive then $\oplus = \max$ and $\otimes = \min$. However we can tackle with this problem partially if we observe that it is sufficient to define $\oplus$ on a domain $D \subseteq [0,1]^2$ such that $(x, y) \in D$ if and only if $x \oplus y \leq 1$ or, by theorems of Acz\'{e}l \cite{aczel}, extend $S$ and $T$ to $[0, \infty)$.

On the other hand observe that in such a logic laws involving sums should be satisfied only if their arguments come from measures of disjoint sets. For example de Morgan laws are too general since they describe a situation of arbitrary sets. For example in a very natural structure $([0,1], +, \cdot, 1 - \mathrm{id})$ which leads to the probability theory de Morgan laws are not satisfied.

Looking for a condition involving $\oplus$ which comes only from disjoint sets we find a following one
\begin{equation*}
(A \cap B) \cup (A \cap B') = A
\end{equation*}
which leads to the definitions

\begin{definition}
A strict negation $n$ is a strictly decreasing function $n : [0,1] \rightarrow [0,1]$ such that $n(0) = 1$ and $n(1) = 0$.
\end{definition}

\begin{definition}
A normal triple $(S, T, n)$ is a triple consisting of a $t$-conorm $S$, $t$-norm $T$ and a strict negation $n$ such that
\begin{equation*}
S(T(x, y), T(x, n(y))) = x
\end{equation*}
for all $x, y \in [0,1]^2$.
We say that a normal triple $(S, T, n)$ is continuous if $S$, $T$ and $n$ are continuous functions.
\end{definition}

\begin{theorem} \label{trojka}
If $(S, T, n)$ is a continuous normal triple, then there exists a unique continuous and strictly increasing function $h : [0,1] \rightarrow [0,1]$ such that $h(0) = 0$, $h(1) = 1$ and
\begin{eqnarray*}
S(x, y) & = & h^{[-1]}(h(x) + h(y)) \\
T(x, y) & = & h^{-1}(h(x) h(y)) \\
n(x) & = & h^{-1}(1 - h(x))
\end{eqnarray*}
In particular $T$ is a strict continuous Archimedean $t$-norm and $S$ is a strict continuous Archimedean $t$-conorm.
\end{theorem}
\begin{proof}
This is theorem 3.2.13 in \cite{alsina}
\end{proof}

\subsection{Measures and fuzzy sets}

Now we postulate the following. Define two fields ($\sigma$-fields) $\mathcal{F}$ and $\mathcal{F}'$ of subsets of $\Omega$ such that if $A \in \mathcal{F}$ then $A' \in \mathcal{F}'$ and vice versa. Sets contained in $\mathcal{F}$ will be called events, while sets in $\mathcal{F}'$ anti-events.

A function $X : (\Omega, \mathcal{F}_\Omega, \mathcal{F}'_\Omega) \rightarrow (E, \mathcal{F}_E, \mathcal{F}'_E)$ is called measurable or a fuzzy variable if for every $A \in \mathcal{F}_E$ and $B \in \mathcal{F}'_E$ we have $X^{-1}(A) \in \mathcal{F}_\Omega$ and $X^{-1}(B) \in \mathcal{F}'_\Omega$.

For a given strict continuous Archimedean $t$-norm $T$ we take a unique $t$-conorm $S$ and a strict negation $n$ such that $(S, T, n)$ constitues a continuous normal triple. For a possibility measure $\Pi : \mathcal{F} \rightarrow [0,1]$ we define a necessity-like measure $\Pi' : \mathcal{F}' \rightarrow [0,1]$ by $\Pi'(A') = n(\Pi(A))$. Then a following diagram commutes

\begin{equation*}
\xymatrix@+20pt{
(\mathcal{F}, \emptyset, X, \subseteq, \cup, \cap) \ar[r]^{\Pi} \ar[d]^{(-)'} & ([0,1], 0, 1, \leq, \max, T) \ar[d]^{n} \\
(\mathcal{F}', \emptyset, X, \subseteq, \cup, \cap) \ar[r]^{\Pi'} & ([0,1], 0, 1, \leq, S, T) \\
}
\end{equation*}

We interpret it as follows. If we know a possibility measure of some event $A$ than we can define a measure of \textit{not} $A$ by $\Pi'(A')$. However if $A' \in \mathcal{F}$ then typically $\Pi(A') \neq \Pi'(A')$.

For $A \in \mathcal{F}$ and $B \in \mathcal{F}'$ sets $A \cup B$ and $A \cap B$ might not belong neither to $\mathcal{F}$ nor $\mathcal{F}'$. However, if $A \cap B = \emptyset$ we can define $\Pi(A \cup B) = \max(\Pi(A), \Pi'(B))$ and similarly if $A$ and $B'$ are $T$-independent, then $\Pi(A \cap B) = T(\Pi(A), \Pi'(B))$.

Finally, some membership functions of fuzzy variables might become degenerated i.e. $\sup A(E)$ may be less than $1$ and should be normalized.
\begin{definition} \label{normalizacja}
Let $T$ be a strict continuous Archimedean $t$-norm with a multiplicative generator $h$. Let $A$ be a degenerated or not fuzzy set in $E$. A normalized version of $A$ denoted by $A^N$ is a fuzzy set defined as
\begin{equation*}
A^N(x) = h^{-1} \left( \frac{1}{\sup h(A(E))} h(A(x)) \right)
\end{equation*}
\end{definition}

In general $\max(A^N, B^N) \neq \max(A, B)^N$. All fuzzy sets should be normalized.

Here we have to mention that this procedure is not well defined and ambiguous. We will show how it works in the proof of Glivenko--Cantelli theorem but we do not know how to tackle with problem of normalization thoroughly.

\section{Theory of fuzzy measurement}

\subsection{Combination of errors}

\begin{proposition}
Let $T$ be a continuous Archimedean $t$-norm and let $A \in \mathcal{FN}(\R)$. Let $(X_n)_{n=1}^\infty$ be a sequence of $T$-independent fuzzy number variables all with distributions given by $A$. Let $Y$ be a fuzzy variable with a distribution given by a characteristic function of some set $B$ and arbitrary dependence with $X_n$ and let $Z_n = X_n + Y$ for all $n \in \N$. Define
\begin{eqnarray*}
A_n & = & Ave(Z_1, \ldots, Z_n) \\
M_n & = & Med(Z_1, \ldots, Z_n)
\end{eqnarray*}
Then $(A_n)_{n=1}^\infty$ and $(M_n)_{n=1}^\infty$ converge in measure and almost surely to $Y$.
\end{proposition}
\begin{proof}
Observe that for any $x_1, \ldots, x_n, y \in \R$ there is 
\begin{eqnarray*}
Ave(x_1 + y, \ldots, x_n + y) & = & Ave(x_1, \ldots, x_n) + y \\
Med(x_1 + y, \ldots, x_n + y) & = & Med(x_1, \ldots, x_n) + y
\end{eqnarray*}
which is also true for intervals. Hence by NFK formula $A_n = \tilde{A}_n + Y$ and $M_n = \tilde{M}_n + Y$ where
\begin{eqnarray*}
\tilde{A}_n & = & Ave(X_1, \ldots, X_n) \\
\tilde{M}_n & = & Med(X_1, \ldots, X_n)
\end{eqnarray*}
By theorem \ref{SLoLN} we know that $\tilde{A}_n$ and $\tilde{M}_n$ converge to $\M A$ in measure. By
\begin{equation*}
\Pi( | A_n - A | \geq \epsilon ) = \Pi( | \tilde{A}_n - \tilde{A} | \geq \epsilon ) \stackrel{n \rightarrow \infty}{\longrightarrow} 0
\end{equation*}
we obtain $A_n \rightarrow A$ in measure. The same holds for $M_n$.
\end{proof}

This proposition allows to treat both random and systematic components of uncertainty in a consistent way in the fuzzy variables theory. A fuzzy number represents a random component and is obtained due to carrying out a measurement. A non-fuzzy set represents a systematic component of uncertainty and is obtained by use of expert method. Such an approach with an arithmetic based on $t$-norms was proposed in \cite{urbanski03} and \cite{urbanski0} and studied empirically in \cite{urbanski08}.

The following section develops a theorem necessary for a construction of empirical membership function of a fuzzy process.

\subsection{Estimation of membership function}

Now we will prove a fuzzy version of Glivenko--Cantelli theorem. Assume we are given a strict continuous Archimedean $t$-norm $T$ with a multiplicative generator $h$. By theorem \ref{trojka} there exists unique $t$-conorm $S$ and a strict negation $n$ such that $(T, S, n)$ is a normal continuous triple. 

For $E = \R^n$ we consider a field of events $\mathcal{F}$ containing all closed subsets of $\R^n$. Similarly a field of anti-events $\mathcal{F}'$ contains all open subsets of $\R^n$. For example if we are given two $T$-independent fuzzy variables $X, Y : \Omega \rightarrow \R$ and two points $x, y \in \R$ then the expression $\Pi(X \leq x, Y < y)$ is equal
\begin{equation*}
\Pi(X \leq x, Y < y) = T(\Pi(X \leq x), \Pi'(Y < y)) = T(\Pi(X \leq x), n(\Pi(Y \geq y)))
\end{equation*}

If $x_1, \ldots, x_n$ is a sequence of real numbers such that $x_1 \leq x_2 \leq \ldots \leq x_n$, then denote $x_{k:n} = x_k$. If $x_1, \ldots, x_n$ is any sequence of real numbers, then $x_{k:n}$ is equal $x_k$ after ordering. Obviously $x_{1:n} = \min(x_1, \ldots, x_n)$ and $x_{n:n} = \max(x_1, \ldots, x_n)$. Under these assumptions and notations we have

\begin{theorem}[Glivenko--Cantelli theorem] \label{histo2n}
Let $T$ be a strict continuous Archimedean $t$-norm with a multiplicative generator $h$ and $(X_n)_{n=1}^\infty$ be a sequence of $T$-independent identically distributed fuzzy number variables with a continuous membership function $A \in \mathcal{FN}(\R)$. Denote $\phi(\alpha) = A^\alpha$ and define estimators of $\phi$ as
\begin{equation*}
\phi_{n}(\alpha) = \left\{ \begin{array}{rl} 
\left[ X_{\left\lceil \frac{n h(\alpha)}{2} \right\rceil : n}, X_{n - \left\lfloor \frac{n h(\alpha)}{2} \right\rfloor : n} \right] & \mathrm{if\;}n\;\mathrm{is\;even} \\
\left[ X_{\left\lceil \frac{(n + 1) h(\alpha)}{2} \right\rceil : n}, X_{n - \left\lfloor \frac{(n + 1) h(\alpha)}{2} \right\rfloor : n} \right] & \mathrm{if\;}n\;\mathrm{is\;odd} \end{array} \right.
\end{equation*}
Then
\begin{enumerate}
	\item $\M \phi_{n}(\alpha) \stackrel{n \rightarrow \infty}{\longrightarrow} \phi(\alpha)$ for every $\alpha \in (0, 1]$,
	\item for every $\epsilon > 0$
\begin{equation*}
D_{n} = \sup_{\alpha \in [\epsilon, 1]} d(\phi_{n}(\alpha), \phi(\alpha))
\end{equation*}
converges to zero in measure and almost surely.
\end{enumerate}
\end{theorem}

\begin{figure}[htb]
	\centering
		\includegraphics[width=0.60\textwidth]{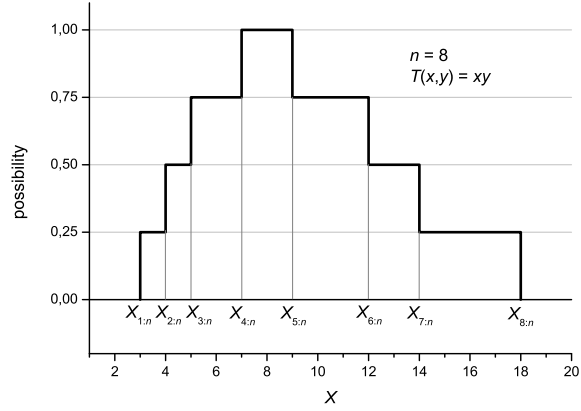}
	\caption{An example of a fuzzy histogram obtained from 8 measurement points.}
	\label{fig:histo}
\end{figure}

The idea of building a membership function of the process is presented in the figure \ref{fig:histo}. Unfortunately, the estimators $\phi_n$ are slightly biased. In order to analyze  this bias we have to prove a more technical version of the theorem.

\begin{theorem}[Glivenko--Cantelli theorem] \label{histo}
Let $T$ be a strict Archimedean $t$-norm with a multiplicative generator $h$ and $(X_n)_{n=1}^\infty$ be a sequence of $T$-independent identically distributed fuzzy number variables with a continuous membership function $A \in \mathcal{FN}(\R)$. Let $b : \N \times (0, 1] \rightarrow \R$ be any function such that
\begin{enumerate}
	\item $\lim_{n \rightarrow \infty} b(n, \alpha) = 0$ for all $\alpha \in (0, 1]$,
	\item $k(n, \alpha) = \frac{n}{2}(h(\alpha) + b(n, \alpha)) \in \{ 1, 2, \ldots, \left\lceil \frac{n}{2} \right\rceil \}$ for all $n \in \N$ and $\alpha \in (0, 1]$,
	\item if $\alpha \leq \beta$ then $k(n, \alpha) \leq k(n, \beta)$ for all $n \in \N$.
\end{enumerate}
Finally denote $\phi(\alpha) = A^\alpha$ and define estimators of $\phi$
\begin{equation*}
\phi_n(\alpha) = [X_{k(n, \alpha) : n}, X_{n - k(n, \alpha) : n}]
\end{equation*}
Then
\begin{enumerate}
	\item $\M \phi_n(\alpha) = \phi \left( h^{-1} ( h(\alpha) + b(n, \alpha) ) \right)$ for every $\alpha \in (0, 1]$,
	\item for every $\epsilon > 0$
\begin{equation*}
D_n = \sup_{\alpha \in [\epsilon, 1]} d(\phi_n(\alpha), \phi(\alpha))
\end{equation*}
converges to zero in measure and almost surely.
\end{enumerate}

\end{theorem}

In this formulation $b$ is a bias which can be added by hand and we see that $\phi_n$ are quite immune to bias. From this theorem we have
\begin{corollary}
Theorem \ref{histo2n} holds with $|b(n, \alpha)| \leq 2/n$.
\end{corollary}

\begin{proof}[Proof of the theorem \ref{histo}]
For a first point set $n \in \N$ and $\alpha \in (0, 1]$ and take any $x = [\underline{x}, \overline{x}] \in \mathcal{I}$. Denote $k = k(n, \alpha)$ and observe that
\begin{eqnarray*}
\Pi( \phi_n(\alpha) = x ) & = & \Pi ( \mathrm{there\;are\;exactly\;}k\;\mathrm{indices\;such\;that\;} X_j \leq \underline{x} \;\mathrm{and\;}k\;\mathrm{such\;that\;}X_j \geq \overline{x} ) \\
& = & \Pi \left( \bigcup_{\sigma \in S_n} \left\{ X_{\sigma(1)} \leq \underline{x}, \ldots, X_{\sigma(k)} \leq \underline{x}, X_{\sigma(n - k + 1)} \geq \overline{x}, \ldots, X_{\sigma(n)} \geq \overline{x}, \right. \right. \\
& & \qquad \left. \left. (X_{\sigma(k + 1)} \leq \underline{x} \vee X_{\sigma(k + 1)} \geq \overline{x})', \ldots, (X_{\sigma(n - k)} \leq \underline{x} \vee X_{\sigma(n - k)} \geq \overline{x})' \right\} \right) \\
& = & \Pi \left( X_1 \leq \underline{x}, \ldots, X_k \leq \underline{x}, X_{n - k + 1} \geq \overline{x}, \ldots, X_n \geq \overline{x} \right. \\
& & \qquad \left.(X_{k+1} \leq \underline{x} \vee X_{k+1} \geq \overline{x})', \ldots, (X_{n-k} \leq \underline{x} \vee X_{n-k} \geq \overline{x})' \right) \\
& = & T \left( \underbrace{\Pi(X_1 \leq \underline{x}), \ldots, \Pi(X_1 \leq \underline{x})}_{k}, \underbrace{\Pi(X_1 \geq \overline{x}), \ldots, \Pi(X_1 \geq \overline{x})}_{k}, \right. \\
&& \left. \underbrace{n(\Pi(X_1 \leq \underline{x} \vee X_1 \geq \overline{x})), \ldots, n(\Pi(X_1 \leq \underline{x} \vee X_1 \geq \overline{x}))}_{n-2k} \right) \\
& = & h^{-1} \left( h^k(A(\underline{x})) \cdot h^k(A(\overline{x})) \cdot \left(1 - h(\max(A(\underline{x}), A(\overline{x}))) \right)^{n-2k} \right)
\end{eqnarray*}

Since $h$ is an increasing homeomorphism we may set $y = h(A(\underline{x}))$ and $z = h(A(\underline{x}))$ and proceed as follows. A degenerated fuzzy set $B_n(y, z) = (y^k z^k (1-\max(y, z))^{n-2k})$ defined on $[0,1]^2$ has a unique maximum at $y_0 = z_0 = \frac{2k}{n}$ for $y, z \in [0, 1]^2$ and by continuity of $A$ and monotonicity of $h$ there exists unique $x = h^{-1}(\phi(2k/n))$ at which maximum is attained. Hence
\begin{equation*}
\M \phi_n(\alpha) = \phi \left( h^{-1} ( h(\alpha) + b(n, \alpha) ) \right) 
\end{equation*}

For a second part by assumption 1 on $b$ and continuity of $h$ and $\phi$ we obtain $\M \phi_n(\alpha) \stackrel{n \rightarrow \infty}{\longrightarrow} A^\alpha$ for all $\alpha \in (0, 1]$.
A normalization of $B_n$ gives
\begin{eqnarray*}
B^N_n(y, z) & = & \frac{n^n}{(2k)^{2k} (n-2k)^{n-2k}} y^k z^k (1-\max(y, z))^{n-2k} \\
& = & \left( \frac{(y z)^{(h(\alpha) + b(n, \alpha))/2} (1-\max(y, z))^{1-h(\alpha)-b(n, \alpha)}}{(h(\alpha) + b(n, \alpha))^{h(\alpha) + b(n, \alpha)}(1-h(\alpha)-b(n, \alpha))^{1-h(\alpha)-b(n, \alpha)}} \right)^n
\end{eqnarray*}
By symmetry we may limit ourselves to the case $z \leq y$. Substituting $y_n = h(\alpha)+b(n, \alpha)$ we find $B^N_n(y_n, y_n) = 1$. Increasing $B^N_n$ by taking $y = z$ we obtain
\begin{equation*}
B^N_n \leq F_n(y) = \left( \frac{y^{h(\alpha) + b(n, \alpha)} (1 - y)^{1 - h(\alpha) - b(n, \alpha)}} {(h(\alpha) + b(n, \alpha))^{h(\alpha) + b(n, \alpha)}(1 - h(\alpha) - b(n, \alpha))^{1 - h(\alpha) - b(n, \alpha)}} \right)^n
\end{equation*}
Taking a derivative of $F_n$ we have
\begin{equation*}
\frac{\partial \ln F_n}{\partial y} = \frac{y - h(\alpha) - b(n, \alpha)}{y(y-1)}
\end{equation*}
so $F_n < 1$ for all $y \in [0, 1] \backslash \{y_n\}$ which means that $B^N_n \stackrel{n \rightarrow \infty}{\longrightarrow} \bs{1}_{\{(h(\alpha), h(\alpha))\}}$ pointwise for all $\alpha \in (0, 1]$. For monotonous convergence take $\tilde{B}^N_n(y, z) = \sup_{k \geq n} B^N_k(y, z)$ and observe that fuzzy sets $B^N_n$ are subsets of $\tilde{B}^N_n$. Since $\tilde{B}^N_n(y, z) \stackrel{n \rightarrow \infty}{\longrightarrow} \bs{1}_{\{(h(\alpha), h(\alpha))\}}$ pointwise, then by a proposition \ref{zbieg} $\phi_n(\alpha) \rightarrow \phi(\alpha)$ in distribution.

For uniform convergence apply the standard procedure. We will focus on the left slope because the procedure for the right one is the same. Take $\epsilon > 0$ and $M \in \N$. Define points
\begin{equation*}
x_{M,k} = \underline{A}^\epsilon + \frac{\overline{A}^\epsilon - \underline{A}^\epsilon}{M} k
\end{equation*}
for $k = 0, 1, \ldots, M$. and $\alpha_k = A(x_k)$. If $\alpha$ is such that $\alpha \in [\alpha_k, \alpha_{k+1}]$ then we have
\begin{equation*}
\underline{\phi}_n(\alpha) - \underline{\phi}(\alpha) \leq \underline{\phi}_n (\alpha_{k+1}) - \underline{\phi} (\alpha_k) \leq \underline{\phi}_n (\alpha_{k+1}) - \underline{\phi} (\alpha_{k+1}) + \frac{1}{M}
\end{equation*}
as well as $\underline{\phi}_n(\alpha) - \underline{\phi}(\alpha) \geq \underline{\phi}_n (\alpha_k) - \underline{\phi} (\alpha_k) - \frac{1}{M}$. Hence we have
\begin{equation*}
| \underline{\phi}_n(\alpha) - \underline{\phi}(\alpha)| \leq \max ( \underline{\phi}_n (\alpha_{k+1}) - \underline{\phi} (\alpha_{k+1}), \underline{\phi}_n (\alpha_k) - \underline{\phi} (\alpha_k) ) + \frac{1}{M} = \max(\underline{\Delta}^{(+)}_{n, M, k} + \underline{\Delta}^{(0)}_{n, M, k}) + \frac{1}{M}
\end{equation*}
where we defined two $\underline{\Delta}$s for simplicity. Thus we have
\begin{equation*}
\underline{D}_n = \max( \max_{0 \leq k \leq M-1} \underline{\Delta}^{(+)}_{n, M, k}, \max_{1 \leq k \leq M} \underline{\Delta}^{(0)}_{n, M, k} ) + \frac{1}{M}
\end{equation*}
We have shown that any $\underline{\Delta}$ converges to zero, hence finite maxima converge to zero. Thus we have $\limsup_{n \rightarrow \infty} \underline{D}_n \leq 1/M$. The same reasoning for $\overline{D}_n$ completes a proof.

\end{proof}

Note that if $F : \R \rightarrow [0,1]$ is a continuous cumulative distribution of a random variable $X$ then a corresponding fuzzy distribution $A \in \mathcal{FN}(\R)$ is given by
\begin{equation*}
A^\alpha = \left[ \sup F^{-1} \left( \frac{\alpha}{2} \right), \inf F^{-1} \left( 1 - \frac{\alpha}{2} \right) \right]
\end{equation*}
This is a probability--possibility transformation similar to those presented in \cite{dubois04} such that $A^\alpha = I_\alpha$ where $I_\alpha$ is a $\alpha$-confidence interval defined as the right hand side of the equation above.

\subsection{Realization of fuzzy variables}

In the theory of probability we can simulate realizations of random variables by means of inverse distribution function. A similar theorem holds for fuzzy variables. 

\begin{proposition} \label{los}
Let $U : \Omega \rightarrow [0, 2]$ be a fuzzy variable with a membership function $A_U(x) = x \bs{1}_{[0, 1]}(x) + (2 - x) \bs{1}_{(1, 2]}(x)$. Let $X$ be a fuzzy number variable with a continuous membership function $A_X \in \mathcal{FN}(\R)$ and define a function $\psi : [0, 2] \rightarrow \R$
\begin{equation*}
\psi(\alpha) = \begin{multi} 
A_X^{-1}(\alpha) & \mathrm{if} & \alpha \leq 1 \\
A_X^{-1}(2 - \alpha) & \mathrm{if} & \alpha > 1
\end{multi}
\end{equation*}
This function is continuous and $\psi(U) = X$ almost surely.
\end{proposition}
\begin{proof}
We have $\Pi(\psi(U) = x) = \Pi(U = \psi^{-1}(x))$ which for $x \leq \M A_X$ is equal $\Pi(U = A(x)) = A(x)$. For $x > \M A_X$ we have $\Pi(U = \psi^{-1}(x)) = \Pi(U = 2 - A(x)) = A(x)$.
\end{proof}

This proposition is a generalization of probability--possibility transformations and shows that only fuzzy numbers represent probabilistic quantities. It allows to simulate fuzzy processes and, together with Glivenko--Cantelli theorem and the Strong Law of Large Numbers shows, that the theory is consistent in this sense, that it allows to produce realizations of fuzzy variables and then regain their distributions by means of membership function estimators.

\section{Conclusions and generalizations}

Theorem \ref{histo} allows to construct a membership function of a fuzzy process from empirical data by a direct estimation of its $\alpha$-cuts. This theorem guarantees that the presented procedure described in the theorem converges to a real membership function.

If empirical data are clearly probabilistic, then a multiplication is a good choice of a $t$-norm. However for real processes another choice of a $t$-norm might be more suitable. The problem of estimation of a $t$-norm was raised and analyzed empirically in \cite{urbanski08}. The procedure is based on minimalization of a distance between fuzzy sets obtained by empirical averaging and theoretical one with use of a $t$-norm. A distance of two fuzzy sets is given by
\begin{equation*}
d(A, B) = \sup_{\alpha \in [\epsilon, 1]} d_{\mathcal{K}(X)}(A^\alpha, B^\alpha)
\end{equation*}
for some set $\epsilon > 0$.

Now we must point out some important problems encountered in the procedure of normalization and proof of Glivenko--Cantelli theorem. The procedure of normalization developed in this paper seems to be ambiguous. The second problem is encountered when a generalization of Glivenko--Cantelli theorem is considered. The assumption about continuity of a membership function of a process is important and cannot be rejected. In order to make a possibility of use non-strict $t$-norms we should weaken a definition of strict negation to the following one.
\begin{definition}
A negation is a function $n : [0,1] \rightarrow [0,1]$ which is continuous, strictly decreasing and such that $n(1) = 0$.
\end{definition}
With this set-up it is possible to generalize theorem \ref{trojka} and obtain all Archimedean $t$-norms. However then the proof of Glivenko--Cantelli theorem breaks down since $h$ is not a homeomorphism.

\end{document}